\documentclass[10pt]{article}

\usepackage{amsmath}
\usepackage{amssymb}
\usepackage{amsthm}
\usepackage{amsrefs}
\usepackage{mathtools}
\usepackage{xcolor}
\usepackage{ytableau}
\usepackage{subcaption}
\usepackage{tikz}
\usepackage{url}

\definecolor{CELL_C}{RGB}{223,223,223}
\definecolor{CELL_B}{RGB}{163, 163, 163}
\newcommand\qbin[3]{\left[ \begin{matrix} #1 \\ #2 \end{matrix} \right]_{\displaystyle #3}}

\newtheorem{theorem}{Theorem}[section]
\newtheorem{proposition}[theorem]{Proposition}
\newtheorem{corollary}[theorem]{Corollary}
\newtheorem*{remark}{Remark}

\author{Hunter Waldron}
\title{Bijective Approaches for Schmidt-Type Theorems}

\date{}
\begin{document}

\maketitle

\begin{abstract}
We provide new Schmidt-type results through an investigation of two bijections, which are
results involving partitions with parts counted only at given indices.
Mork's bijection, the first of these, was originally
given as a proof of Schmidt's theorem. We show that a version of Sylvester's bijection is
equivalent to Mork's bijection applied to 2-modular diagrams, which implies refinements of existing
results and new generating function identities. 
We then develop a bijection 
based on an idea appearing in a recent paper of Andrews and Keith, that places partitions
counted at the indices $r$, $t+r$, $2t+r, \dots$ in correspondence with $t$-colored partitions.
This leads to a substantial generalization of an identity of Bridges and Uncu, and complements a
similar investigation of Li and Yee.
\end{abstract}

\section{Introduction}

A \textit{partition} is a sequence $\lambda = (\lambda_1, \lambda_2, \lambda_3, \dots)$ of weakly decreasing
non-negative integers that eventually become zero. The \textit{parts} of $\lambda$ are the non-zero terms,
the number of which is written $\ell(\lambda)$. When convenient, we write $\lambda = (\lambda_1, \dots,
\lambda_{\ell(\lambda)})$, ignoring the trailing zeroes. The \textit{size} of $\lambda$, written $\lvert \lambda \rvert$,
is the sum of all  parts, $\lambda_1 + \lambda_2 + \lambda_3 + \cdots$. If $\lambda$ has size $n$,
we say $\lambda$ is a partition of $n$. For example, the 7 partitions are 5 are

\begin{figure}[!ht]
\centering
(5), (4,1), (3,2), (3,1,1), (2,2,1), (2,1,1,1), and (1,1,1,1,1).
\end{figure}

\noindent The first 3 of these partitions are the partitions of 5 with \textit{distinct parts}, meaning no part appears more
than once.

There has been rapidly growing interest \cite{andrews_keith, li_yee, bridges_uncu, berkovich_uncu, ji, andrews_paule, uncu, chern_yee}
in the research of partitions counted with \textit{Schmidt-weights}, where  parts are only counted at specified indices.
These are so named after Frank Schmidt who in \cite{schmidt} made the remarkable observation:

\begin{theorem}[Schmidt \cite{schmidt}]
The number of partitions of $n$ is equal to the number of partitions $\lambda$ with distinct parts such that
$\lambda_1 + \lambda_3 + \lambda_5 + \cdots = n$.
\end{theorem}

For example,
the 7 partitions with distinct parts corresponding to the partitions listed above are

\begin{figure}[!ht]
\centering
(5), (5,1), (5,2),  (5,3), (5,4), (4, 2, 1), and (4, 3, 1).
\end{figure}

Relationships between families of partitions like Schmidt's often translate into identities for generating functions.
Let $\mathcal{P}$ be the set of all partitions and $\mathcal{D}$ be the set of partitions with distinct parts. Then
Schmidt's theorem is equivalent to the equality
$$
\sum_{\lambda \in \mathcal{D}} q^{\displaystyle \lambda_1 + \lambda_3 +
\lambda_5 + \cdots} = \sum_{\lambda \in \mathcal{P}} q^{\displaystyle \lvert \lambda \rvert}.
$$
The sum on the right is well-known to be equal to $1/(q;q)_\infty$, giving the identity
$$ \sum_{\lambda \in \mathcal{D}} q^{\displaystyle \lambda_1 + \lambda_3 +
\lambda_5 + \cdots} = \frac{1}{(q;q)_\infty}. $$
This is written using the standard notation for the q-Pochhammer symbol, shown below.
When possible, we use standard definitions and notation such as this throughout the paper. More can be found in section 2.
\begin{equation*}
(z; q)_n = \prod_{k=0}^{n-1} (1-zq^k) \hspace{0.25cm} \text{for integers $n \geq 0$}, \hspace{0.5cm} (z; q)_\infty = \prod_{k=0}^\infty (1-zq^k).
\end{equation*}

One of the many known proofs of Schmidt's theorem is a bijection given by Mork \cite{mork},
which matches certain hooks from a given partition to the distinct parts of another. Beyond
an application by Chern and Yee in \cite{chern_yee}, this map has not been studied in much
depth.

We show in section 3 that
a well-known version of Sylvester's bijection by Bessenrodt \cite{bessenrodt}, which proves a
refinement of Euler's classic theorem on the equinumerosity of partitions of equal size with
distinct and odd parts, is essentially equivalent to Mork's bijection. Indeed, Mork's bijection
becomes identical when applied to 2-modular diagrams of partitions with odd parts,
forgetting the numbers.
We use the equivalence of these maps to add a Schmidt-type refinement to the already long
list of properties that Bessenrodt's bijection is known to possess.

\begin{proposition} Let $\lambda$ be a partition of $n$ with distinct parts and first part $k$, such that
$\lambda_1 + \lambda_3 + \lambda_5 + \cdots = m$, and let $\mu$ be the preimage of $\lambda$ under
Bessenrodt's bijection. Then $\mu$ has $2m - n$ parts and first part $1 + 2k + 2n - 4m$.
\label{proposition:euler_refinement}
\end{proposition}

\begin{remark} \normalfont This result shows that there is a strong linear relationship between
the Schmidt-weight $\lambda_1 + \lambda_3 + \lambda_5 + \cdots$ and $\ell(\mu)$ for partitions
of fixed size, which is illustrated in table \ref{table:bessenrodt} below.
\end{remark}

\begin{table}[h!]
\begin{center}
\begin{tabular}{ |c|c|c| }
\hline
$\lambda_1 + \lambda_3 + \lambda_5 + \cdots $ & $\lambda$ with distinct parts & $\mu$ with odd parts \\
\hline
7 &  (7) & (1,1,1,1,1,1,1) \\
\hline
6 & (6,1) & (3,1,1,1,1) \\
\hline
5 & (5,2) & (5,1,1) \\
\hline
5 & (4, 2, 1) & (3,3,1) \\
\hline
4 & (4,3) & (7) \\
\hline
\end{tabular}
\end{center}
\caption{The partitions of 7 with distinct and odd parts.}
\label{table:bessenrodt}
\end{table}

As a byproduct of the work leading to Proposition \ref{proposition:euler_refinement}, we also produce a
refinement of Schmidt's theorem:

\begin{corollary}
The number of partitions of $n$ with $\ell$ parts is equal to the number of partitions $\lambda$ with distinct parts
such that $\lambda_1 + \lambda_3 + \lambda_5 + \cdots = n$ and $\lvert \lambda \rvert = 2n - \ell$.
\label{corollary:schmidt_refinement}
\end{corollary}

The following generating function identities are also implied:

\begin{theorem}
\begin{align}
  \sum_{\lambda \in \mathcal{D}} &z^{\displaystyle \lvert\lambda\rvert}
  q^{\displaystyle \lambda_1 + \lambda_3 + \lambda_5
  + \cdots} = \frac{1}{(qz;qz^2)_\infty}
  \label{theorem:mork_product_identity_1} \\
  \sum_{\lambda \in \mathcal{D}}& z^{\displaystyle \lvert\lambda\rvert}
  q^{\displaystyle \lambda_2 + \lambda_4 + \lambda_6
  + \cdots} = \frac{1}{(z;qz^2)_\infty}
  \label{theorem:mork_product_identity_2}
\end{align}
\label{theorem:mork_product_identity}
\end{theorem}

\begin{remark} \normalfont If we factor (\ref{theorem:mork_product_identity_1}) as
$\sum_{n=0}^\infty p_n(z) q^n$, $p_n(z)$ will be a polynomial, and indeed a $z$-analog of the partition function.
\begin{equation}
\frac{1}{(qz;qz^2)_\infty} = 1 + zq + (z^2+z^3)q^2 + (z^3 + z^4 + z^5)q^3 + \cdots
\end{equation}
These polynomials have two interpretations based on this work,

\begin{equation}
p_n(z) \hspace{0.5cm} = \hspace{0.5cm} \smashoperator{\sum_{\substack{\lambda \in \mathcal{D}, \\
\lambda_1 + \lambda_3 + \lambda_5 + \cdots = n}}} z^{\displaystyle \lvert \lambda \rvert}
\hspace{0.5cm} = \hspace{0.5cm}
\smashoperator{\sum_{\substack{\lambda \in \mathcal{P}, \\ \lvert \lambda \rvert = n}}}
z^{\displaystyle 2 n - \ell(\lambda)}.
\end{equation}
The sum on the right is simply a (shifted and reflected) distribution of the partitions of $n$ by number of parts,
giving a connection between the Schmidt-weight and length similar to what we have with
Proposition \ref{proposition:euler_refinement}.
This is also the algebraic equivalent of Corollary \ref{corollary:schmidt_refinement}.

For (\ref{theorem:mork_product_identity_2}), the inclusion of $z$ is necessary for the
coefficients to be finite, since for example any partition with one part will have
Schmidt-weight zero.
\end{remark}

Our final application of Mork's bijection uses a different approach, counting only
some of the hook lengths used in order to give generating function
identities for partitions $\lambda$ with distinct  parts $\lambda_i$
counted when $i \equiv 1 \pmod{4}$. This technique is probably not difficult to
generalize.

\begin{theorem}
\begin{align}
  \smashoperator{\sum_{\substack{\lambda \in \mathcal{D}, \\ \ell(\lambda) \equiv
  0 \; \mathrm{or} \; 3 \\ (\mathrm{mod} \; 4)}}}
& z^{\displaystyle \lambda_1}
  q^{\displaystyle \lambda_1 + \lambda_5 + \lambda_9 + \cdots} =
  1 + \sum_{n=1}^\infty \frac{q^{2n+1 \choose 2} z^{4n-1}}
  {(qz;q)_n^4}
  \label{theorem:mork_sum_identity_1} \\
  \smashoperator{\sum_{\substack{\lambda \in \mathcal{D}, \\ \ell(\lambda) \equiv
  1 \; \mathrm{or} \; 2 \\ (\mathrm{mod} \; 4)}}}
& z^{\displaystyle \lambda_1}
  q^{\displaystyle \lambda_1 + \lambda_5 + \lambda_9 + \cdots} =
  \sum_{n=1}^\infty \frac{q^{2n \choose 2} z^{4n-3}}
  {(qz;q)_n^2 (qz;q)_{n-1}^2}
  \label{theorem:mork_sum_identity_2}
\end{align}
\label{theorem:mork_sum_identity}
\end{theorem}

In section 4, we change our focus to investigating Schmidt-type results for unrestricted partitions. Several
such results have already been found. For instance, in \cite{bridges_uncu}
Bridges and Uncu gave the following two generating function identites.

\begin{theorem}[Bridges-Uncu \cite{bridges_uncu}]
\begin{align}
\sum_{\lambda \in \mathcal{P}} z^{\displaystyle \lambda_1 }q^{\displaystyle \lambda_1 + \lambda_3 +
\lambda_5 + \cdots} &=  \frac{1}{(qz;q)_\infty^2} \label{theorem:bridges_uncu_1} \\
\sum_{\lambda \in \mathcal{P}} z^{\displaystyle \lambda_1 }q^{\displaystyle \lambda_2 + \lambda_4 +
\lambda_6 + \cdots} &= \frac{1}{(1-z)(qz;q)_\infty^2}
\end{align}
\label{theorem:bridges_uncu}
\end{theorem}

Uncu had found the $z=1$ specialization of (\ref{theorem:bridges_uncu_1}) previously in \cite{uncu},
which was also independently found by Andrews and Paule in \cite{andrews_paule}. They added a
combinatorial interpretation in terms of 2-colored partitions, which Ji \cite{ji} later proved, and
indeed refined, bijectively.

A more general connection with colored partitions was more recently given by Li and Yee in \cite{li_yee}, which they proved
bijectively as well. We have translated their result into the language of colored partitions, which they call
$t$-multipartitions. From here, let $t$ and $r$ be arbitrary positive integers.

\begin{theorem}[Li-Yee \cite{li_yee}]
The number of partitions $\lambda$ such that $\ell(\lambda) = (s-1)t + j$ and
$\lambda_1 + \lambda_{t+1} + \lambda_{2t+1} + \cdots = n$ is equal to the number of $t$-colored
partitions of $n$ where $s$ is the most number of times any color appears, and $j$ is the largest such color.
\label{theorem:li_yee}
\end{theorem}

By developing an idea that was remarked on by Andrews and Keith in \cite{andrews_keith}, we give a bijection,
which we call the \textit{color-conjugate} map, that implies a similar result.

\begin{theorem}
The number of partitions $\lambda$ such that $\lambda_1 = k_1$, $\lambda_r = k_r$,
$\lambda_r + \lambda_{t+r} + \lambda_{2t+r} + \cdots = n$, and for each $1 \leq i \leq t$, 
$$ c_i(\lambda) = \sum_{k \geq 0} (\lambda_{kt + r + i - 1} - \lambda_{kt + r + i}), $$
is equal to the number of pairs $(\nu, \mu)$ such that $\nu$ is a partition with first part
$k_1 - k_r$ and at most $r-1$ parts, and $\mu$ is a $t$-colored partition with $k_r$ parts
where each color $i$ appears $c_i(\lambda)$ times.
 \label{theorem:color_conjugate}
\end{theorem}

\begin{remark} \normalfont
Both of these theorems imply that the partitions $\lambda$ such that $\lambda_1 + \lambda_{t+1}
+ \lambda_{2t+1} + \cdots = n$ are equinumerous with the $t$-colored partitions of $n$, however
 the refinements are distinct. Indeed, when $r=1$, it can
be easily checked that the color-conjugate map places partitions with $(s-1)t + j$ parts in bijection with
$t$-colored partitions where, instead of the condition in Li and Yee's theorem,
$s$ is the largest part size, and $j$ is the largest color used on the largest part size.

For example, when $t=3$, the color-conjugate bijection sends
$(4,4,3,3,3,3)$ to the 3-colored partition $(2^3,2^3,2^3,1^2)$, while the map given by Li and Yee
sends this same partition to $(3^3,3^3, 1^2)$.
\end{remark}

There are two notable interpretations of Theorem \ref{theorem:color_conjugate} as generating
function identities.
Now, Theorem \ref{theorem:bridges_uncu} becomes the special case of
(\ref{theorem:color_conjugate_identity_1}) with $t=2$ and $r \in \{1,2\}$.

\begin{theorem}
\begin{align}
\sum_{\lambda \in \mathcal{P}} z^{\displaystyle \lambda_1} q^{\displaystyle \lambda_r +
\lambda_{t+r} + \lambda_{2t+r} + \cdots} &= \frac{1}{(1-z)^{r-1} (qz;q)_\infty^t}
\label{theorem:color_conjugate_identity_1} \\
\sum_{\lambda \in \mathcal{P}}  \prod_{i=1}^{t}  z_i^{\displaystyle c_i(\lambda)} \times
q^{\displaystyle \lambda_1 + \lambda_{t+1} + \lambda_{2t+1} + \cdots} &=
\prod_{i=1}^{t} \frac{1}{(qz_i; q)_\infty}
\label{theorem:color_conjugate_identity_2}
\end{align}
\label{theorem:color_conjugate_identity}
\end{theorem}

We can actually improve (\ref{theorem:color_conjugate_identity_1}) by adding a
third variable $s$ that tracks the ordinary size of the partitions summed over, if
we abandon our purely bijective approach. This gives us our main result:

\begin{theorem}
\begin{equation}
\sum_{\lambda \in \mathcal{P}} z^{\displaystyle \lambda_1}
s^{\displaystyle \lvert \lambda \rvert} q^{\displaystyle \lambda_r + \lambda_{t+r}
+ \lambda_{2t+r} + \cdots} = 
\frac{1}{(sz;s)_{r-1}} \prod_{n=0}^\infty
\frac{1}{(s^{nt+r} q^{n+1} z; s)_t}
\end{equation}
\label{theorem:main}
\end{theorem}

Replacing $s$ with $q$ and $q$ with $q^{-1}$  in Theorem \ref{theorem:main},
then simplifying the right hand side gives us another interesting generating function identity,
with a straight forward combinatorial interpretation in terms of $2$-colored partitions.

\begin{corollary}
\begin{equation}
\sum_{\lambda \in \mathcal{P}} z^{\displaystyle \lambda_1}
   q^{\displaystyle |\lambda| -  \lambda_r -
   \lambda_{t+r} - \lambda_{2t+r} - \cdots}
   = \frac{1}{(qz;q)_\infty(q^{r-1}z; q^{t-1})_\infty}.
\end{equation}
\end{corollary}

\begin{corollary}
If $r > 1$, the number of partitions with first part $k$ such that all parts except at the
indices $r, t+r, 2t+r, \dots$ sum to $n$ is equal to the number of $2$-colored partitions
of $n$ with $k$ parts where the second color only appears on parts of size
$r-1, t-1 + r-1, 2(t-1) + r-1, \dots$
\label{corollary:opposite_schmidt}
\end{corollary}

\begin{remark} \normalfont
We can also prove Corollary \ref{corollary:opposite_schmidt} by modifying color-conjugate map
such that the parts $\lambda_r, \lambda_{t+r}, \lambda_{2t+r}, \dots$
are used to color the parts $\lambda_{r-1}$, $\lambda_{t-1 + r-1}$, $\lambda_{2(t-1) + r-1}$, \dots.
\end{remark}

\section*{Acknowledgments}

The author would like to thank his advisor William Keith for introducing him to this area of research
and for critiquing early versions of this work.

\section{Preliminaries}

Here we provide additional background and the notation that will be
used in the proofs of section 3 and 4.

For any positive integer $t$, a \textit{$t$-colored partition} is a
partition $\lambda$ with each part $\lambda_i$ marked by a \textit{color}
$c_i \in \{1, \dots, t\}$, with the order of colors across parts of
equal size weakly decreasing. That is, if $1 \leq i < j$ and
$\lambda_i = \lambda_j$, then $c_i \geq c_j$. The notation we use for
$t$-colored partitions places the color as a superscript:
$\lambda = (\lambda_1^{c_1}, \lambda_2^{c_2}, \lambda_3^{c_3}, \dots)$.
The \textit{size} of $\lambda$ is simply $\lambda_1 + \lambda_2
+ \lambda_3 + \cdots$, ignoring the colors.
For example, the 2-colored partitions of 3 are
\begin{equation*}
\begin{gathered}
(3^1), \; (3^2), \; (2^1, 1^1), \; (2^2, 1^1), \; (2^1, 1^2), \;
(2^2, 1^2), \; (1^1, 1^1, 1^1), \\
(1^2, 1^1, 1^1), \; (1^2, 1^2, 1^1), \; \textrm{and} \; (1^2, 1^2, 1^2).
\end{gathered}
\end{equation*}

The \textit{Young diagram} of a partition $\lambda$ is a visual
representation of $\lambda$ as a grid of squares
arranged into rows that are aligned on the left, with the number of squares
in the $j$th row down corresponding to the value of $\lambda_j$. The
\textit{cell} $(i, j)$ of $\lambda$ is the $i$th square down
in the $j$th column from the left in the Young diagram of $\lambda$.
The cells $(i,i)$ form the \textit{diagonal} of $\lambda$.
Swapping $i$ and $j$ in all the cells of $\lambda$ produces another
partitions $\lambda'$, the \textit{conjugate} of $\lambda$.

The \textit{Durfee square} of $\lambda$ is the largest square of cells
that can be placed inside the Young diagram of $\lambda$. Its
\textit{size} is the number of cells in the width, or equivalently in the height.

To any cell $(i,j)$ in the Young diagram of $\lambda$ we can associate a
\textit{hook} which consists of the cell $(i,j)$ and any below and to the right.
The \textit{hook length} is the number of cells in the hook.

\begin{figure}[!ht]
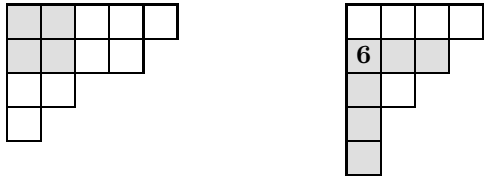

\centering
\ytableausetup{boxsize=1.25em}
\begin{ytableau}
*(CELL_C) & *(CELL_C) &  & & \\
*(CELL_C) & *(CELL_C) &  & \\
& \\ \;
\end{ytableau}
\hspace{2cm}
\begin{ytableau}
\; & & & \\
*(CELL_C)\textbf{6} & *(CELL_C) & *(CELL_C) \\
*(CELL_C) & \\
*(CELL_C) \\
*(CELL_C)
\end{ytableau}
\caption{The partitions $(5,4,2,1)$ and $(4,3,2,1,1)$ are conjugates. The Durfee square is
shown in $(5,4,2,1)$, and the hook of $(2,1)$ is shown in $(4,3,2,1,1).$}
\end{figure}

\begin{figure}[!ht]
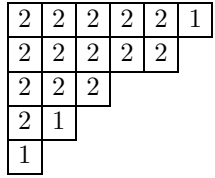

\centering
\ytableausetup{boxsize=1.25em}
\begin{ytableau}
2 & 2 & 2 & 2 & 2 & 1 \\ 2 & 2 & 2 & 2 & 2 \\ 2 & 2 & 2 \\ 2 & 1 \\ 1
\end{ytableau}
\caption{The 2-modular diagram of the partition $(11, 10, 6, 3, 1).$}
\end{figure}

For any integer $m \geq 2$, if we write the parts of a partition $\lambda$ as
$\lambda_i = m \nu_i + r_i$ where $1 \leq r_i \leq m$ and
$\nu = (\nu_1, \nu_2, \nu_3, \dots)$ is a partition, the \textit{m-modular} diagram of $\lambda$
is the Young diagram of $\nu$ with the cell at the end
of the $i$th row filled with $r_i$, and all others with $m$.

The $q$-binomial coefficient is defined using the $q$-Pochhammer symbol
for all integers $n, k \geq 0$ as
\begin{equation*}
\qbin{n}{k}{q} = \frac{(q;q)_n}{(q;q)_k (q;q)_{n-k}},
\end{equation*}
which is useful for giving us an identity for partitions with
Young diagrams whose cells fit inside an $n \times k$ rectangle.
\begin{equation*}
\sum_{\substack{\lambda \in \mathcal{P}, \\ \lambda_1 \leq n, \;
\mathrm{and} \; \\ \ell(\lambda) \leq k}}
\!\!\!\!\!\!\!q^{\displaystyle |\lambda|}= \qbin{n+k}{k}{q}.
\end{equation*}

\section{Mork's bijection}

Mork's bijection \cite{mork} $\mathrm{M} \colon \mathcal{P} \to \mathcal{D}$ is defined on the partition $\mu$ as follows.
For $i \geq 1$, $\mathrm{M}(\mu)_{2i - 1}$ is the $i$th diagonal hook length of $\mu$, and $\mathrm{M}(\mu)_{2i}$
is the hook length of the cell to the right of that. Since the diagonal hook lengths of $\mu$ sum to $\lvert \mu \rvert$,
the odd indexed parts of $\mathrm{M}(\mu)$ sum to $\lvert \mu \rvert$, which implies Schmidt's theorem.
This is easily understood by example.

\begin{figure}[!ht]
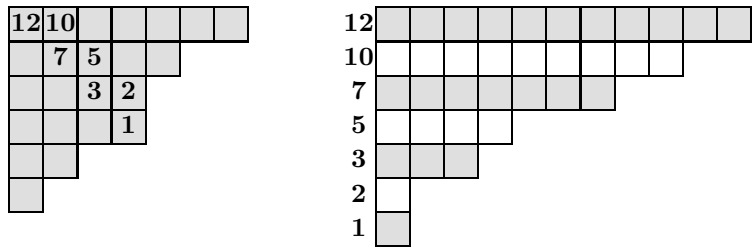

\centering
\ytableausetup{boxsize=1.25em}
\begin{ytableau}[*(CELL_C)]
\textbf{12} & \textbf{10} & & & & & \\ & \textbf{7} & \textbf{5} & & \\
& & \textbf{3} & \textbf{2} \\ & & & \textbf{1} \\ & \\ \;
\end{ytableau}
\hspace{1cm}
\begin{ytableau}
\none[\textbf{12}] & *(CELL_C) & *(CELL_C) & *(CELL_C) & *(CELL_C) & *(CELL_C) &
*(CELL_C) & *(CELL_C) & *(CELL_C) & *(CELL_C) & *(CELL_C) & *(CELL_C) \\
\none[\textbf{10}] & & & & & & & & & \\
\none[\textbf{7}] & *(CELL_C) & *(CELL_C) & *(CELL_C) & *(CELL_C) & *(CELL_C) &
*(CELL_C) & *(CELL_C) \\
\none[\textbf{5}] & & & & \\
\none[\textbf{3}] & *(CELL_C) & *(CELL_C) & *(CELL_C) \\
\none[\textbf{2}] & \\
\none[\textbf{1}] & *(CELL_C)
\end{ytableau}
\caption{Mork's bijection maps the partition $(7,5,4,4,2,1) \in \mathcal{P}$
to $(12,10,7,5,3,2,1) \in \mathcal{D}$. }
\label{figure:mork}
\end{figure}

By similar reasoning, the even indexed parts of $\mathrm{M}(\mu)$ sum to $\lvert \mu \rvert - \ell(\mu)$, so
$\vert \mathrm{M}(\mu) \rvert = 2 \lvert \mu \rvert - \ell(\mu)$. Therefore if we fill the cells of $\mu$
at the end of each row with a 1, and all others with 2s, we construct the 2-modular diagram of
a partition with odd parts that has size equal to $\mathrm{M}(\mu)$. Since this construction defines a
bijection $\mathrm{D} \colon \mathcal{P} \to \mathcal{O}$, $\mathrm{M} \circ \mathrm{D^{-1}}$ is a size preserving bijection
from $\mathcal{O}$ to $\mathcal{D}$. This also immediately implies Corollary \ref{corollary:schmidt_refinement}.

Bessenrodt's bijection, given in \cite{bessenrodt} as a version of Sylvester's bijection, similarly acts on a 2-modular
diagram's diagonal hooks to give a size preserving bijection $\mathrm{B} \colon \mathcal{O} \to \mathcal{D}$. 
For $\mu \in \mathcal{O}$ and $i \geq 1$, Bessenrodt's map assigns to $\mathrm{B}(\mu)_{2i-1}$ the number of 2s and 1s
in the $i$th diagonal hook of $\mu$'s 2-modular diagram, and to $\mathrm{B}(\mu)_{2i}$ the number of 1s.
Although defined differently, we can easily see that $\mathrm{B} = \mathrm{M} \circ \mathrm{D}^{-1}$.
For, from their respective definitions, we immediately see that
$\mathrm{B}(\mu)$ and $\mathrm{M} \circ \mathrm{D}^{-1} (\mu) $ coincide at the odd indices.

As for even indices $2i$, the number of 2s in the cells at and below $(i,i)$ in the 2-modular diagram of $\mu$
must be equal to the number of cells at and below $(i,i+1)$
since any 2 must have a cell to the right, but a 1 may not. Similarly the cells to the
right of $(i,i)$ with 2s will be one less than the number of all such cells.
Thus both maps coincide at the even indices.
We use this equivalence to prove the following. \newline

\textbf{Proof of Proposition \ref{proposition:euler_refinement}.} Let $\mu^* = \mathrm{D}^{-1}(\mu)$. Then Mork's bijection maps $\mu^*$ to $\lambda$. So,
as above, $2 \lvert \mu^* \rvert - \ell(\mu^*) = \lvert \lambda \rvert$. Since $\lvert \mu^* \rvert = m$ and
$\ell(\mu^*) = \ell(\mu)$ the first property follows.

For the second property, clearly the first diagonal hook length of $\mu^*$ is $\lambda_1$, and
$\mu^*_1 = (\mu_1 - 1)/2$, so

\begin{equation}
\begin{aligned}
\lambda_1 = \mu^*_1 + \ell(\mu^*) - 1 &= \frac{\mu_1 - 1}2 + \ell(\mu) - 1 \\
&= \frac{\mu_1 - 1}2 + 2m - n - 1.
\end{aligned}
\end{equation}
Solving for $\mu_1$ completes the proof.

\qed

\subsection*{Identities implied by Mork's bijection}

\textbf{Proof of Theorem \ref{theorem:mork_product_identity}.}
From the work above, Mork's bijection implies the equality

\begin{equation}
\sum_{\lambda \in \mathcal{D}} z^{\displaystyle \lvert \lambda \rvert} q^{\displaystyle \lambda_1 +
\lambda_3 + \lambda_5 + \cdots} = \sum_{\mu \in \mathcal{P}} z^{\displaystyle 2 \lvert \mu \rvert - \ell(\mu)}
q^{\displaystyle \lvert \mu \rvert}.
\end{equation}

The sum on the right has a simple interpretation. For any partition $\mu \in \mathcal{P}$ with $n$ parts,
the first $n$ cells on left in $\mu$'s Young diagram are counted as $(qz)^n$, and to the right of that may
be any partition with at most $n$ parts, with each cell contributing $qz^2$. Summing over $n$ gives

\begin{equation}
\sum_{\mu \in \mathcal{P}} z^{\displaystyle 2 \lvert \mu \rvert - \ell(\mu)}
q^{\displaystyle \lvert \mu \rvert} = \sum_{n=0}^\infty \frac{(qz)^n}{(qz^2; qz^2)_n}
= \sum_{n=0}^\infty \frac{(qz^2)^n}{(qz^2; qz^2)_n} (z^{-1})^n.
\end{equation}

Using the well-known identity
$$ \sum_{n=0}^\infty \frac{q^n}{(q; q)_n} z^n = \frac{1}{(qz;q)_\infty} $$
with $q$ and $z$ replaced with $qz^2$ and $z^{-1}$, respectively, this becomes
\begin{equation}
\frac{1}{(z^{-1} qz^2; qz^2)_\infty}=  \frac{1}{(qz; qz^2)_\infty}
\end{equation}
proving (\ref{theorem:mork_product_identity_1}).

To prove (\ref{theorem:mork_product_identity_2}), we can use a similar argument, starting instead with.
\begin{equation}
\sum_{\lambda \in \mathcal{D}} z^{\displaystyle \lvert \lambda \rvert} q^{\displaystyle \lambda_2 +
\lambda_4 + \lambda_6 + \cdots} = \sum_{\mu \in \mathcal{P}} z^{\displaystyle 2 \lvert \mu \rvert - \ell(\mu)}
q^{\displaystyle \lvert \mu \rvert - \ell(\mu)}.
\end{equation}

\qed

\textbf{Proof of Theorem \ref{theorem:mork_sum_identity}.} From the definition of Mork's bijection,
if $\lambda \in \mathcal{D}$, then $\ell(\lambda) \equiv 0 \; \textrm{or} \; 3 \pmod{4}$ exactly when
the Durfee square of $\mathrm{M}^{-1}(\lambda)$ is has even length. Moreover, we have that
the parts $\lambda_1, \lambda_5, \lambda_9, \dots$ are given by the 1st, 3rd, 5th, \dots diagonal hook
lengths of $\mathrm{M}^{-1}(\lambda)$, which here we will denote $h_1^\mu$, $h_3^\mu$,
$h_5^\mu, \dots$. Using these facts, the following equality must hold.
\begin{equation}
 \smashoperator{\sum_{\substack{\lambda \in \mathcal{D}, \\ \ell(\lambda) \equiv
  0 \; \mathrm{or} \; 3 \\ (\mathrm{mod} \; 4)}}}
 z^{\displaystyle \lambda_1}
  q^{\displaystyle \lambda_1 + \lambda_5 + \lambda_9 + \cdots} =
 \smashoperator{\sum_{\substack{\mu \in \mathcal{P}\; \textrm{with} \\ \textrm{even Durfee length}}}}
 z^{\displaystyle h_1^\mu} q^{\displaystyle h_1^\mu + h_3^\mu + h_5^\mu + \cdots}
\end{equation}

Consider a partition $\mu$ with Durfee length $2n$. The cells that are both in $\mu$'s Durfee square and in
the hooks counted by the above sum are in total the $2n$th triangular number, contributing
$q^{2n+1 \choose 2}z^{4n-1}$. This can be seen by deleting the other cells, then sliding these hooks up
and to the left.

To the right of $\mu$'s Durfee square and, in conjugate, below,
may be any partition with at most $2n$ parts where only the odd indexed parts contribute to $q$, since these are the
remainder of the hook lengths, and the first part contributes to $z$. The identity
\begin{equation}
\smashoperator{\sum_{\substack{\lambda \in \mathcal{P}, \\ \ell(\lambda) \leq 2n}}} z^{\displaystyle \lambda_1}
q^{\displaystyle \lambda_1 + \lambda_3 + \lambda_5 + \cdots} = \frac{1}{(qz;q)_n^2}
\end{equation}
can easily be derived from the color-conjugate bijection in section 4. This all together implies that

\begin{equation}
\smashoperator{\sum_{\substack{\mu \in \mathcal{P}\; \textrm{with} \\ \textrm{even Durfee length}}}}
 z^{\displaystyle h_1^\mu} q^{\displaystyle h_1^\mu + h_3^\mu + h_5^\mu + \cdots} =
1 + \sum_{n=1}^\infty \frac{q^{2n+1 \choose 2}z^{4n-1}}{(qz;q)_n^4}
\end{equation}
which gives (\ref{theorem:mork_sum_identity_1}). The proof of (\ref{theorem:mork_sum_identity_2}) is similar.

\qed

\section{The color-conjugate map}

Let $\mathcal{P}^t$ be the set of all $t$-colored partitions, and $\mathcal{P}_{\leq r-1}$ be the set of partitions with at most
$r-1$ parts.
We define the color-conjugate map $\mathrm{C}_{t,r} \colon \mathcal{P} \to  \mathcal{P}_{\leq r - 1} \times \mathcal{P}^t$ 
as follows for given $t, r \geq 1$ on a partition $\lambda$.

Let $\nu$ be the partition given by $\nu_i = \lambda_i - \lambda_r$ for $1 \leq i \leq r$. The
$t$-colored partition $\mu$ is constructed as follows. The parts of $\mu$ are given by the conjugate of the partition
$(\lambda_r, \lambda_{t+r}, \lambda_{2t+r}, \dots)$. The color of $\mu_i$ is given by subtracting
$r-1$ from $\lambda$'s $i$th column, then counting the number of cells from the bottom of the column up, stopping
at the first counted part. We then define $\mathrm{C}_{t,r}(\lambda) = (\nu, \mu)$.
\begin{figure}[!ht]
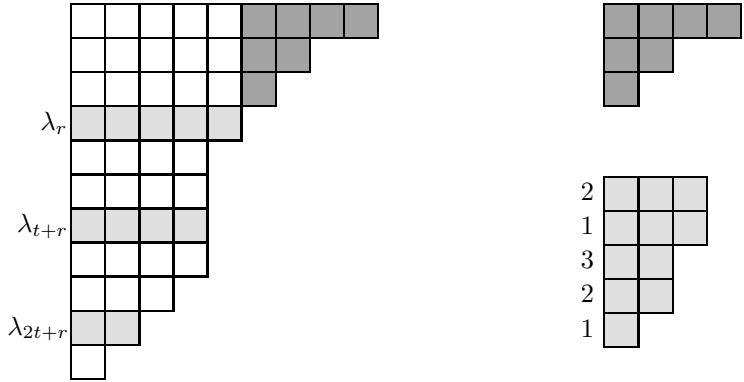

\ytableausetup{boxsize=1.25em}
\begin{subfigure}[t]{0.49\textwidth}
\centering
$$\begin{ytableau}
\none & & & & & &*(CELL_B) & *(CELL_B) & *(CELL_B) & *(CELL_B) \\
\none & & & & & &*(CELL_B) & *(CELL_B) \\
\none & & & & & & *(CELL_B) \\
\none[\lambda_r] & *(CELL_C) & *(CELL_C) & *(CELL_C) & *(CELL_C) & *(CELL_C) \\
\none & & & & \\
\none & & & & \\
\none[\!\!\!\!\! \lambda_{t+r}] & *(CELL_C) & *(CELL_C) & *(CELL_C) & *(CELL_C) \\
\none & & & & \\
\none & & & \\
\none[\!\!\!\!\!\!\! \lambda_{2t+r}]  & *(CELL_C) & *(CELL_C) \\
\none &
\end{ytableau} $$
\end{subfigure}
\hfill
\begin{subfigure}[t]{0.49\textwidth}
\centering
$$ \begin{ytableau}[*(CELL_B)]
\; &  &  & \\
\; & \\
\;
\end{ytableau} $$
\vspace{0,5cm}
$$ \begin{ytableau}
\none[2] & *(CELL_C) & *(CELL_C) & *(CELL_C) & \none & \none \\
\none[1] & *(CELL_C) & *(CELL_C) & *(CELL_C) \\
\none[3] & *(CELL_C) & *(CELL_C) \\
\none[2] & *(CELL_C) & *(CELL_C) \\
\none[1] & *(CELL_C)
\end{ytableau} $$
\end{subfigure}
\caption{$\mathrm{C}_{3,4}$ maps the partition $\lambda = (9,7,6,5,4,4,4,4,3,2,1)$ to $\nu=(4,2,1)$ and $\mu=(3^2, 3^1, 2^3, 2^2, 1^1)$.}
\end{figure}

The claimed properties immediately follow from the construction.
We can give $\mathcal{C}_{t,r}^{-1}$ just as easily, showing that $\mathcal{C}_{t,r}$ is indeed a bijection.

Given $\nu$ and $\mu$, let $\lambda_i = \ell(\mu) + \nu_ i$ for $1 \leq i \leq r-1$. For $1 \leq i \leq \ell(\mu)$,
append $(t-1)\mu_i + c_i$ to $\lambda$'s $i$th column where $c_i$ is the color of $\mu$'s $i$th part.
\newline

\textbf{Proof of Theorem \ref{theorem:color_conjugate_identity}.}
Using the well-known identities
\begin{equation}
 \smashoperator{\sum_{\lambda \in \mathcal{P}_{\leq r-1}}} z^{\displaystyle \lambda_1} = \frac{1}{(1-z)^{r-1}}  \hspace{0.5cm} \text{and}
 \hspace{0.5cm} \sum_{\lambda \in \mathcal{P}^t}  z^{\displaystyle \lambda_1}  q^{\displaystyle \ell(\lambda)} = \frac{1}{(qz;q)_\infty^t},
\end{equation}
(\ref{theorem:color_conjugate_identity_1}) is  immediately implied from the properties of the color-conjugate map, and
(\ref{theorem:color_conjugate_identity_2}) similarly.

\qed

\subsection*{Proof of the main result}

We first prove Theorem \ref{theorem:main} in the $r=1$ case. That result will then be used to
find the general case. For integers $n \geq 0$ define
\begin{equation}
   f_n(q,s) = \sum_{\substack{\lambda \in \mathcal{P}, \\ \lambda_1 = n}}
   s^{\displaystyle |\lambda|}
   q^{\displaystyle \lambda_1 + \lambda_{t+1} + \lambda_{2t+1} + \cdots}.
\end{equation}
The first step of the proof is finding a recurrence relation for $f_n(q, s)$.

Consider any partition $\lambda$ with first part $n$. The first part will
be counted as $(qs)^n$ by $f_n(q,s)$. For each possibility
$0 \leq k = \lambda_{t+1} \leq n$, in $\lambda$'s Young diagram
we are guaranteed a $k \times (t-1)$
rectangle of cells between $\lambda_1$ and $\lambda_{t+1}$ that is counted as
$s^{k(t-1)}$, and an $(n-k) \times (t-1)$ rectangle to the right
that may contain any partition. This is counted as
$$ \qbin{n - k + t - 1}{t - 1}{s} $$
and finally, the part $\lambda_{t+1}$ and beyond is counted by $f_k(q,s)$.
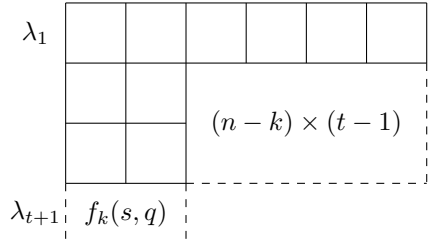
\begin{figure}[!ht]
\centering
\begin{tikzpicture}[scale=0.8]
\draw (0,0) -- (6,0);
\draw (0,-2) -- (2,-2);
\draw (0,-3) -- (2,-3);
\draw (0,-1) -- (6,-1);
\draw (0,0) -- (0,-3);
\draw (1,0) -- (1,-3);
\draw (2,0) -- (2,-3);
\draw (3,0) -- (3,-1);
\draw (4,0) -- (4,-1);
\draw (5,0) -- (5,-1);
\draw (6,0) -- (6,-1);
\draw[dashed] (0,-3) -- (0, -4);
\draw[dashed] (2,-3) -- (2, -4);
\draw[dashed] (2,-3) -- (6,-3);
\draw[dashed] (6,-1) -- (6,-3);
\node at (1,-3.5) {$f_k(s, q)$};
\node at (4,-2) {$(n-k) \times (t-1)$};
\node at (-0.5,-0.5) {$\lambda_1$};
\node at (-0.5,-3.5) {$\lambda_{t+1}$};
\end{tikzpicture}
\caption{A visualization of the cells counted by $f_n(q,s)$. In this
example $n=6$, $t=3$, and $k=2$.}
\end{figure}

Some elementary combinatorial reasoning then gives us the recurrence
\begin{equation}
   f_n(q,s) = (qs)^n \sum_{k=0}^n s^{k(t-1)} \qbin{n - k + t - 1}{t - 1}{s}
   f_k(q,s).
\end{equation}
Now let $F(s, q, z)$ be the generating function for $f_n(q, s)$:
$$ F(s, q, z) = \sum_{n=0}^\infty f_n(q, s) z^n. $$
This is indeed the sum appearing on the left
in the statement of Theorem \ref{theorem:main} since the power of
$z$ coincides with the first part of the partitions counted.

Replacing $f_n(q,s)$ with the recurrence found above in $F(s, q, z)$,
\begin{equation}
\begin{split} F(s, q, z) &=
   \sum_{n=0}^\infty \left((qs)^n \sum_{k=0}^n s^{k(t-1)}
   \qbin{n - k + t - 1}{t - 1}{s} f_k(q,s) \right)z^n  \\
   &= \sum_{n=0}^\infty \sum_{k=0}^n s^{k(t-1)} 
   \qbin{n - k + t - 1}{t - 1}{s} f_k(q,s) (sqz)^n.
\end{split}
\end{equation}
Recognizing this as the Cauchy product of two power series in the variable
$sqz$, we may factor the above to get
\begin{equation}
F(s,q,z) = \left( \sum_{n=0}^\infty s^{n(t-1)} f_n(q, s) (sqz)^n \right)
   \left( \sum_{n=0}^\infty \qbin{n + t - 1}{t - 1}{s} (sqz)^n \right).
\end{equation}
The factor on the left is plainly $F(s, q, s^t qz)$ after combining
powers. The factor on the right may be
interpreted as an instance of the generating function
\begin{equation}
\sum_{\substack{\lambda \in \mathcal{P}, \\ \ell(\lambda) \leq t}}\!
(qz)^{\displaystyle \lambda_1} s^{\displaystyle |\lambda|}.
\end{equation}

A simpler identity for this is $1/(sqz; s)_t$. Thus
\begin{equation}
 F(s, q, z) = \frac{1}{(sqz; s)_t} F(s, q, s^tqz)
\end{equation}
which is a functional equation that can easily be solved for $F(s, q, z)$.
Applying this relation $N$ times produces
\begin{equation}
   F(s,q, z) = \prod_{n=0}^N \frac{1}{(s^{nt+1} q^{n+1} z;s)_t}
   F(s,q,s^{Nt} q^N z).
\end{equation}
Letting $N \to \infty$ scales the terms of $F(s, q, s^{Nt} q^N z)$
by increasingly large powers of $s^t q$. In the limit, we obtain
\begin{equation}
F(s,q, z) = \prod_{n=0}^\infty \frac{1}{(s^{nt+1} q^{n+1} z;s)_t}
\end{equation}
which proves the $r=1$ case.

From here, assume $r > 1$ and define
\begin{equation}
   F_r(s, q, z) = \sum_{\lambda \in \mathcal{P}} z^{\displaystyle \lambda_1}
   s^{\displaystyle |\lambda|} q^{\displaystyle \lambda_r +
   \lambda_{t+r} + \lambda_{2t+r} + \cdots}.
\end{equation}

For any integer $n \geq 0$ consider the partitions
$\lambda$ with $\lambda_r = n$.
The $n \times (r - 1)$ rectangle of cells
above $\lambda_r$ is counted by $F_r(s, q, z)$ as $s^{n(r-1)} z^n$ and
the $\infty \times (r-1)$ rectangle to the right of that may contain
any partition with at most $r-1$ parts, counted as
\begin{equation}
\sum_{\substack{\lambda \in \mathcal{P}, \\ \ell(\lambda) \leq r-1}}\!\!\!\!
\!\!z^{\displaystyle \lambda_1} s^{\displaystyle |\lambda|},
\end{equation}
which is simply $1/(sz; s)_{r-1} $.
The parts $\lambda_r$ and below are counted by $f_n(s, q)$.
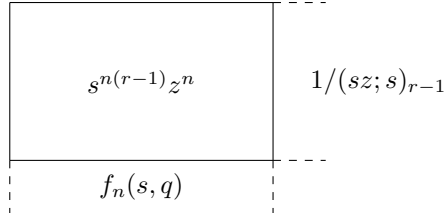
\begin{figure}[!ht]
\centering
\begin{tikzpicture}[scale=0.7]
\draw (0,0) -- (5, 0);
\draw (0,-3) -- (5, -3);
\draw (0,0) -- (0, -3);
\draw (5,0) -- (5, -3);
\draw[dashed] (5, 0) -- (6, 0);
\draw[dashed] (5, -3) -- (6, -3);
\draw[dashed] (0,-3) -- (0, -4);
\draw[dashed] (5,-3) -- (5, -4);
\node at (2.5, -1.5) {$s^{n(r-1)}z^n $};
\node at (2.5, -3.5) {$f_n(s, q)$};
\node at (7, -1.5) {$ 1/(sz;s)_{r-1}$};
\end{tikzpicture}
\caption{A visualization of what is counted by $F_r(s,q,z)$ for
a given $n$.}
\end{figure}

Summing over $n$ produces 
\begin{equation}
\begin{split}
    F_r(s, q, z) &= \sum_{n=0}^\infty s^{n(r-1)} z^n
   \frac{1}{(sz;s)_{r-1}} f_n(s, q) =
   \frac{1}{(sz;s)_{r-1}} F(s, q, s^{r-1} z) \\
   &= \frac{1}{(sz;s)_{r-1}} \prod_{n=0}^\infty
   \frac{1}{(s^{nt+r} q^{n+1} z; s)_t}
\end{split}
\end{equation}
which is the desired result.
\qed

\section*{Further work}

The research of Schmidt-weighted partitions is still a very new idea in partition theory, and new results
and techniques are appearing rapidly. Much of the low hanging fruit has already been picked, but there are
still many unexplored ideas.

As for what we considered in this paper,
generalizing Mork's bijection in any way would be interesting. Maybe the most obvious approach is to first generalize Bessenrodt's map,
perhaps to $m$-modular diagrams with multiple allowed remainders. The natural way to this while preserving size
is to count the number of 1s, 2s, 3s, \dots, $m$s in the first diagonal hook, then 2s, 3s, \dots, $m$s, then
3s, \dots, $m$s, and so on, giving the parts of the resulting partition, and repeating for each diagonal hook.
This will not in general be a bijection unless the partitions considered are restricted in some way.
For example, the following partitions written as 3-modular diagrams both map to (5,4,3,1):

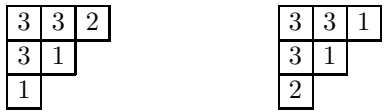
\begin{figure}[!ht]
\centering
\ytableausetup{boxsize=1.25em}
\begin{ytableau}
3 & 3 & 2 \\ 3 & 1 \\ 1
\end{ytableau}
\hspace{2cm}
\begin{ytableau}
3 & 3 & 1 \\ 3 & 1 \\ 2
\end{ytableau}
\end{figure}

The color-conjugate map can be used to describe partitions with restrictions on individual colors, since
on the Schmidt-weight side these become the alternating sums in Theorem \ref{theorem:color_conjugate}.
This can be used to find results similar to those considered in \cite{andrews_keith}. One could also attempt to modify
the proof of Theorem \ref{theorem:main} to work with different Schmidt-weights or to sum over a different set
of partitions.

\bibliography{final_draft}

\begin{bibdiv}
\begin{biblist}

\bib{andrews_keith}{article}{
      author={Andrews, George~E.},
      author={Keith., William~J.},
       title={Schmidt-type theorems for partitions with uncounted parts},
        date={2022},
        note={\url{https://arxiv.org/abs/2203.05202}},
}

\bib{andrews_paule}{article}{
      author={Andrews, George~E.},
      author={Paule, Peter},
       title={Macmahon's partition analysis xiii: Schmidt type partitions and
  modular forms},
        date={2022},
     journal={Journal of Number Theory},
      volume={234},
       pages={95\ndash 119},
}

\bib{berkovich_uncu}{article}{
      author={Berkovich, Alexander},
      author={Uncu, Ali},
       title={On finite analogs of schmidt's problem and its variants},
        date={2022},
        note={\url{https://arxiv.org/abs/2205.00527}},
}

\bib{bessenrodt}{article}{
      author={Bessenrodt, Christine},
       title={A bijection for lebesgue's partition identity in the spirit of
  sylvester},
        date={1994},
     journal={Discrete Mathematics},
      volume={132},
       pages={1\ndash 10},
}

\bib{bridges_uncu}{article}{
      author={Bridges, Walter},
      author={Uncu, Ali},
       title={Weighted cylindric partitions},
        date={2022},
     journal={Journal of Algebraic Combinatorics},
}

\bib{chern_yee}{article}{
      author={Chern, Shane},
      author={Yee, Ae~Ja},
       title={Diagonal hooks and a schmidt-type partition identity},
        date={2022},
     journal={The Electronic Journal of Combinatorics},
      volume={29},
      number={P2.10},
}

\bib{ji}{article}{
      author={Ji, Kathy~Q.},
       title={A combinatorial proof of a schmidt type theorem of andrews and
  paule},
        date={2022},
     journal={The Electronic Journal of Combinatorics},
      volume={29},
}

\bib{li_yee}{article}{
      author={Li, Runqiao},
      author={Yee, Ae~Ja},
       title={Schmidt type partitions},
        date={2022},
     journal={Enumerative Combinatorics and Applications},
      volume={3},
}

\bib{mork}{article}{
      author={Mork, Peter},
       title={Interrupted partitions, solution to problem 10629},
        date={2000},
     journal={The American Mathematical Monthly},
      volume={107},
       pages={974},
}

\bib{schmidt}{article}{
      author={Schmidt, Frank},
       title={Problem 10629: Interrupted partitions},
        date={1999},
     journal={The American Mathematical Monthly},
      volume={104},
       pages={87\ndash 88},
}

\bib{uncu}{article}{
      author={Uncu, Ali},
       title={Weighted rogers–ramanujan partitions and dyson crank},
        date={2017},
     journal={Ramanujan Journal},
      volume={46},
       pages={579 \ndash  591},
}

\end{biblist}
\end{bibdiv}

\end{document}